\documentclass[12pt]{article}
\usepackage{color,amsmath,amsthm,amssymb}
\usepackage{hyperref}

\newtheorem{thm}{Theorem}

\newtheorem{pro}{Proposition}

\theoremstyle{definition}
\newtheorem{const}{Construction}

\newcommand{\ol}{\overline}

\newcommand{\mbb}{\mathbb}
\newcommand{\mc}{\mathcal}

\newcommand{\vs}{\vskip 15pt}
\newcommand{\vhs}{\vskip 6pt}

\title{Group-Case Commutative Association Schemes and Their Character Tables}
\author{Sung Y. Song \\ {\small Department of Mathematics, Iowa State University, Ames, Iowa, 50011, U. S. A.} \\[5mm] Hajime Tanaka \\ {\small Graduate School of Information Sciences, Tohoku University, Sendai, 980-8579, Japan}}
\date{}

\begin{document} 

\maketitle

\begin{abstract}
\noindent
Leading towards the classification of primitive commutative association schemes as the ultimate goal, Bannai and some of his school have been trying to
\begin{itemize}
\item identify the major sources of (primitive) commutative association schemes,
\item collect known group-case primitive commutative association schemes, and
\item compute their character tables
\end{itemize}
over the last twenty years.
The construction of their character tables are important first step for a systematic study of such association schemes and towards the classification of those schemes.
In this talk, we briefly survey the progress made in this direction of research, and list some open problems.
\end{abstract}
\vs

\begin{center}{\bf{Dedicated to Eiichi Bannai on the occasion of his 60th
birthday}}\end{center}

\section{Introduction}

Let a finite group $G$ act on a finite set $X$ transitively.
Then $G$ naturally acts on $X\times X$ by $(x, y)^g = (x^g, y^g)$.
Let $R_0, R_1, \dots, R_d$ be the orbits of $G$ on $X\times X$ with $R_0
= \{(x, x) \ :\ x\in X\}$.
Then $\mathcal{X}=(X, \{R_i\}_{0\le i\le d})$ is an association scheme, called a {\it group-case} (or {\it Schurian}) association scheme, and denoted by $\mathcal{X}(G, X)$.\vhs

Let $G$ be a finite group, and let $H$ be a subgroup of $G$. Let
$X=H\backslash G$ be the set of the cosets of $H$ in $G$. Then $G$
acts transitively on $X$ under the action $(Hx)^g=H(xg)$. The
group-case scheme $\mc{X}(G, H\backslash G)$ is commutative if and
only if the permutation character of $G$ on $H\backslash G$ is
multiplicity-free. Any group-case scheme $\mc{X}(G,X)$ can be viewed
as $\mc{X}(G, H\backslash G)$ with the point stabilizer $H=G_x$ for
an $x\in X$. The condition that the group-case association scheme is
primitive is equivalent to that of the permutation group acts on the
cosets by a maximal subgroup.\vhs

In early nineteen eighties, Bannai had a conviction that the works
by many group theorists on the classification of maximal subgroups
of finite simple groups (by using the classification of finite
simple groups) would eventually lead to the complete list of
group-case primitive commutative association schemes. He seemed to
believe that the calculations of parameters and character tables of
known association schemes of such kind were to be feasible. He began to
investigate the major sources of group-case primitive commutative
association schemes, collect examples, and calculate their character
tables.\vhs

The major sources of group-case (primitive) commutative association
schemes of large class which Bannai (cf. \cite{Bi90}) has considered
were as follows.

\begin{enumerate}
\item The actions of classical groups (or Chevalley groups) $G$ on appropriate
subspaces $X$ of vector spaces over finite fields. (If the subspaces
are isotropic, they are well understood because then the groups act
on the cosets by parabolic subgroups. The cases of non-isotropic
subspaces require further study. See examples in Table 1 below.)

\item The actions of classical groups (or Chevalley groups) $G$ on the
cosets $X=H\backslash G$ of multiplicity-free (maximal) subgroups
$H$. (A pair $(G,H)$ of a finite group $G$ and a subgroup $H$ whose
permutation character $1_H^G$ is multiplicity-free, is often called
a Gelfand pair. See examples in Tables 2 and 3 below.)

\item Finite (simple) groups, loops, and quasigroups $G$. (Commutative
association schemes obtained from these algebraic structures by
using their conjugacy classes as we will see in the sequel.)

\item The $n$-dimensional vector spaces $V$ over $GF(q)$ and
subgroups $H$ of $GL(n,q)$; for example, the action of the
semidirect product $G=V\rtimes H$ on $V$, in other words, $G$ acts on
$X=H\backslash G$.
\end{enumerate}\vhs

We note that the items in this list are not necessarily mutually
exclusive nor cover all primitive commutative association schemes.
Character tables of many commutative association schemes coming from
the permutation groups in the above list have been investigated by
many people including,
Bannai-Song \cite{BS89a,BS89b,BS90},
Bannai-Shen-Song \cite{BSS90},
Bannai-Kawanaka-Song \cite{BKS90},
Bannai-Kwok-Song \cite{BKwS90},
Bannai-Shen-Song-Wei \cite{BSS91},
Kwok \cite{Kw91, Kw92},
Henderson \cite{He01a,He01b,He03,He01c},
Tanaka \cite{Ta01,Ta02},
Bannai-Tanaka \cite{BT03},
Fujisaki \cite{Fu03,Fu06},
and Bannai-Song-Yamada \cite{BSY06}.\vhs

In Section \ref{sec:tables}, we briefly recall the definition of character tables of
commutative association schemes and related basic facts. In Section
\ref{sec:loops}, we discuss Paige's simple Moufang loops and their character
tables. In Section \ref{sec:group-case}, we give a list of known examples of group-case
association schemes whose character tables are either calculated or
conjectured. Some open cases coming from Gelfand pairs are listed.
In Section \ref{sec:G}, we illustrate another example that is constructed in a
different way from the previous ones. The character tables
illustrated in this note are happened to be closely related to each
other in an interesting way.\vhs

Our aim is to review the progress that has been made in the construction
of character tables by paying attention to the construction methods
employed. In doing this we would like to point out some connections
between association schemes and classical groups and geometries.
Experts in algebraic combinatorics, groups and geometry will
hopefully find some useful information and sufficient pointers in
this note.

\section{The character tables}\label{sec:tables}
Let $\mathcal{X}=(X, \{R_i\}_{0\le
i\le d})$ be a commutative association scheme. Let $A_0, A_1, \dots,
A_d$ be the adjacency matrices and let $E_0, E_1, \dots, E_d$ be the
primitive idempotents of $\mathcal{X}$. The Bose-Mesner algebra
$\mathcal{A}=\langle A_0, A_1, \dots, A_d \rangle = \langle E_0,
E_1, \dots, E_d \rangle$ of $\mathcal{X}$ over the field $\mathbb{C}$
of complex numbers, satisfies
\[A_j = \sum\limits_{i=0}^{d}p_j(i)E_i.\]

\noindent Equivalently, with the character table $P=\left
[p_j(i)\right ]$ of $\mathcal{X}$, we have
\[\left [ A_0\ A_1\ \cdots\ A_d\right ]\ =\ \left [E_0\ E_1\ \cdots\ E_d \right ] \left [ \begin{array}{ccccc}
1 & k_1 &k_2 &\cdots &  k_d\\
1 & p_1(1)&p_2(1) &\cdots & p_d(1)\\
1 & p_1(2)&p_2(2) &\cdots & p_d(2)\\
\vdots & \vdots & \vdots & \ddots &\vdots \\
1& p_1(d) &p_2(d) &\cdots & p_d(d)\\
\end{array}\right ]
.\]

\noindent The character table $P$ of the association scheme
satisfies the (i) row and (ii) column orthogonality relations
\cite[Theorem 3.5]{BI84}. For $ i,j\in \{0,1,\dots, d\},$
\[\mbox{(i)}\quad \sum\limits_{l=0}^d
\frac{1}{k_l}p_l(i)\ol{p_l(j)}=\frac{|X|}{m_i}\delta_{ij},\qquad
\mbox{(ii)}\quad \sum\limits_{l=0}^d
m_lp_i(l)\ol{p_j(l)}=|X|k_i\delta_{ij},
\] where $m_l=\mbox{ rank}( E_l) =\mbox{ trace}(E_l)\ (0\le l\le
d)$, $\delta_{ij}$ is the Kronecker delta, and $\ol{a}$ denotes the
complex conjugate of $a$. The numbers $m_l$ are the multiplicities
of the scheme.

\vs Let $G$ be a finite (simple) group, and let $C_0, C_1, \dots,
C_d$ be the conjugacy classes of $G$. Then by defining the associate
classes $R_i$ by
\[(x, y)\in R_i \quad {\mbox { iff }}\quad yx^{-1}\in C_i,\] we
obtain a (primitive) commutative scheme $\mathcal{X}(G)=(G,
\{R_i\}_{0\le i\le d})$ which is often referred to as the {\it group
scheme}. Given a finite (simple) quasigroup, we also obtain a
(primitive) commutative association scheme by defining the associate
relations as above. For $i=0, 1, \dots, d$, let $k_i$ and $f_i$
denote the sizes of conjugacy classes $C_i$ and the degrees of the
irreducible characters of $G$, respectively. Then
$f_i=\sqrt{\mbox{rank}(E_i)}$ and the character table $P$ of
$\mathcal{X}(G)$ and the group character table $T$ has the following
relation:
\[T=\left [ \begin{array}{cccc} f_0 & &  & 0\\ & f_1  & & \\
&  & \ddots & \\ 0 & & & f_d\\ \end{array} \right ]\cdot P\cdot
\left [ \begin{array}{cccc} 1/k_0 & &  & 0\\ & 1/k_1  & & \\
&  & \ddots & \\ 0 & & & 1/k_d\\ \end{array} \right ].
\]

\section{Character tables of Paige's Moufang loops}\label{sec:loops}

A set $Q$ with one binary operation is a \textit{quasigroup} if the
equation $xy=z$ has a unique solution in $Q$ whenever two of
$x,y,z\in Q$ are specified. A \textit{loop} is a quasigroup with a
neutral element $1$ satisfying $1x=x=x1$ for every $x\in Q$.

A \textit{Moufang loop} is a loop in which any of the following
(equivalent) Moufang identities holds: $((xy)x)z=x(y(xz))$,\quad
$x(y(zy))=((xy)z)y$,\quad $(xy)(zx)=x((yz)x)$,  or
$(xy)(zx)=(x(yz))x$.

Paige (1956) \cite{Pa56} introduced a class of finite simple Moufang
loops which we are referring to as Paige's simple Moufang loops.
Liebeck (1987) \cite{Li87} proved that there are no other finite
(non-associative) simple Moufang loops besides these Moufang loops.
For every finite field $\mbb{F}_q$, there is exactly one simple
Moufang loop $\mc{M}^*=\mc{M}^*(q)$ of order $q^3(q^4-1)/(q-1,2)$.
Bannai and Song (1989) \cite{BS89b} calculated the character tables
of $\mc{M}^*(q)$. We now recall the definition of $\mc{M}^*(q)$.

On the set
\[\left \{\left [\begin{array}{cc} a & \alpha \\ \beta & b\\
\end{array}\right ]:\  a,b\in \mbb{F}_q,\ \alpha, \beta \in
\mbb{F}_q^3 \right \}\]
where $\mbb{F}_q^3:=\{(x,y,z):\ x,y,z\in \mbb{F}_q\},$ using the dot
product $\alpha\cdot\beta$ and vector product $\alpha\times\beta$ in
$\mbb{F}_q^3$, {{define the Zorn's multiplication
\[ \left [\begin{array}{cc}
a & \alpha \\ \beta & b\\ \end{array}\right ]
\left [\begin{array}{cc} c & \gamma \\ \delta & d\\
\end{array}\right ]:=\left [\begin{array}{cc} ac+\alpha\cdot\delta &
a\gamma +d\alpha - \beta\times\delta
\\ c\beta+b\delta+\alpha\times\gamma
 & \beta\cdot\gamma+ bd\\
\end{array}\right ].\]}}

\noindent Given $M= \left [\begin{array}{cc} a & \alpha \\ \beta & b\\
\end{array}\right ]$, we define its determinant by \[\mbox{det}(M) :=
ab-\alpha\cdot\beta.\] Then both sets $\mc{L}:=\{M:\mbox{
det}(M)\neq 0\}$ and $\mc{M}:=\{M:\mbox{ det}(M)=1\}$ are
(non-associative) Moufang loops.
The center of $\mc{L}$ is \[Z(\mc{L})=\left \{\left
[\begin{array}{cc} a& 0\\ 0 & a\\ \end{array}\right ]\ :\ a\in
\mbb{F}_q^*\right \}.\] We see that $|Z(\mc{M})|=1$ if the
characteristic of $\mbb{F}_q$ is 2, otherwise, $|Z(\mc{M})|=2$. The
quotient loop $\mc{M}^*:=\mc{M}/Z(\mc{M})$ is referred to as the
\textit{Paige's simple Moufang loop.}

\medskip
The conjugacy classes and character tables of $\mc{M}^*(q)$ were
calculated in \cite{BS89b}. Here we recall the case with $q=2^r$.
\begin{thm}{\rm{\cite[Table 5]{BS89b}}}
The character table of $\mc{X}(\mc{M^*})$, $q=2^r$ is given by
\[ \left [ \begin{array}{cc|ccc|ccc}
1 &(q^6-1) &q^6-q^3& \cdots &q^6-q^3 &q^6+q^3 &\cdots &q^6+q^3\\
1 & q^2-1 & -q^3+q^2 & \cdots & -q^3+q^2 & q^3+q^2 & \cdots & q^3+q^2\\
\hline 1 &-q^3-1 & & & & & &\\
\vdots & \vdots &  & q^2\left [a_{kl}\right ]& & & 0 & \\
 1 &-q^3-1 & & & & & &\\
\hline  1 &q^3-1 & & & & & &\\
\vdots &\vdots & & 0 & & & q^2\left [b_{mn}\right ] &\\
 1 &q^3-1 & & & & & &\\
\end{array} \right ]\]
where
\[a_{kl}=-q(\sigma^{kl}+\sigma^{-kl}), \ 1\le k,l\le q/2;\ \sigma=\exp(2\pi i/(q+1))\]
\[b_{mn}=q(\rho^{mn}+\rho^{-mn}),\ 1\le m,n\le (q-2)/2;\ \rho=\exp(2\pi i/(q-1)).\]
\end{thm}
\vhs

\noindent This character table resembles the following table of
$\mc{X}(PSL(2,q))$, $q=2^r$ which is derived from the group
character table of $PSL(2,2^r)$ found in \cite{Do71}.
\begin{pro}{\rm{\cite[\S 38]{Do71}}} The character table of $\mc{X}(PSL(2,q))$,
$q=2^r$ is given by
\[ \left [ \begin{array}{cc|ccc|ccc}
1 &(q^2-1) &q^2-q& \cdots &q^2-q &q^2+q &\cdots &q^2+q\\
1 & 0 & -q+1 & \cdots & -q+1 & q+1 & \cdots & q+1\\
\hline 1 &-q-1 & & & & & &\\
\vdots & \vdots &  & \left [a_{kl}\right ]& & & 0 & \\
 1 &-q-1 & & & & & &\\
\hline  1 &q-1 & & & & & &\\
\vdots &\vdots & & 0 & & &\left [b_{mn}\right ] &\\
 1 &q-1 & & & & & &\\
\end{array} \right ]\]
where $a_{kl}$ and $b_{mn}$ as in the above theorem.
\end{pro}\vhs
So, it is evident how the character table of $\mc{X}(\mc{M^*})$ can
be expressed in terms of $\mc{X}(PSL(2,q))$, and vice versa. It is
also shown that the corresponding result for odd $q$ is very similar
to the above. Namely, if we replace all $q^3$ by $q$ and $q^2$ by 1
in the character table of $\mc{X}(\mc{M}(q)^*)$ for $q=p^r$ with $p$
an odd prime, then the resulting table is the character table of the
fusion scheme $\tilde{\mc{X}}(PSL(2,q))$ obtained from that of
$\mc{X}(PSL(2,q))$ by combining two conjugacy classes of order $p$
of $PSL(2,q)$, which are conjugate in $PGL(2,q)$, into a single
class. (See \cite[Theorem 2.3.3]{BS89b}.)

\begin{const}
The character table of
$\mc{X}(\mc{M}(q)^*)$ was constructed as follows.\\
(1) We calculated the parameters $p_{ij}^h$ of $\mc{X}(\mc{M}(q)^*)$
and expressed them in terms of parameters
of $\mc{X}(PSL(2,q))$ as in \cite[Lemma 2.2.2]{BS89b}.\\
(2) Using the relationship
$\sum\limits_{h=0}^dp_{ij}^hp_h(r)=p_i(r)p_j(r)$ between the
parameters and the characters for $\mc{X}(PSL(2,q))$, and the
relationship between the two sets of parameters obtained in (1),
we found the relations between the parameters and entries of the
character table of $\mc{X}(\mc{M}(q)^*)$.\\
(3) Finally, we examined if the table satisfied the orthogonality
conditions of rows and columns to be a character table.
\end{const}
Of course, this method only works when two schemes are intimately
related so that the relations in (1) are simple enough to figure out
the relations in (2). The construction of the character tables of
Paige's Moufang loops has led us to be able to construct those for
many other association schemes via this method. Despite the
`intimacy' requirement, it has been used effectively in many cases
where one scheme `controls' many others. See, for example
\cite{BSS90, BSS91, Ta01}.

\section{Group-case association schemes}\label{sec:group-case}
All association schemes coming from the major sources listed in
Introduction are essentially coming from either a finite simple
groups or quasigroups, Gelfand pairs, the primitive permutation
groups and the suitable subgroups that can be the point stabilizers
of the actions of groups. All group schemes $\mc{X}(G)$ can be
viewed as groups $G\times G$ acting on $G$ by $x\mapsto g^{-1}xh$.
For given a quasigroup $Q$ and $x\in Q$, Suppose $G$ is the
multiplicative group $Gr(Q)$ of $Q$ generated by all permutations
$L(x)$ and $R(x)$ of $Q$ defined by
\[L(x):\ y\mapsto xy; \quad R(x):\ y\mapsto yx.\]
Then the association scheme defined by the orbits of $Gr(Q)$ on
$Q\times Q$ as the associate relations, is isomorphic to the
quasigroup association scheme $\mc{X}(Q)$ defined by its conjugacy
classes.

The class of group schemes and quasigroup schemes contain all finite
simple groups and quasigroups. All finite simple groups have been
classified (cf. \cite{GL94}), but quasigroups seems to require a lot
more work. Many maximal subgroups of finite simple groups have been
discovered by the work of many group theorists (cf., \cite{In88,
ILS86}). So there are a lot of group-case primitive commutative
association schemes to be studied. The character tables of some of
these association schemes have been calculated. In doing so, the
following facts play an important role.

\begin{const}\label{const:B}
Let $G$ be a group and $H$ be a
subgroup of $G$. Let $c_0,c_1,\dots,c_d$ be class representatives of
the conjugacy classes $C_0,C_1,\dots, C_d$. Let $\{Hg_{j}H :\ 0\le j\le d\}$
be the set of all double cosets of $H$.
Suppose \[1_H^G = \rho_0 + \rho_1 + \cdots + \rho_d\] is
the decomposition of $1_H^G$ into irreducible characters $\rho_0,
\rho_1, \dots, \rho_d$ of $G$. Then by \cite[Corollary 11.7]{BI84},
the entries $p_j(i)$ of the character table of $\mc{X}(G,H\backslash G)$
are given by
\begin{align*}
p_{j'}(i) & =\frac{1}{|H|}\sum\limits_{c\in Hg_{j}H}\rho_{i}(c)\\
 & = \frac{1}{|H|}\sum\limits_{k}|Hg_{j}H\cap C_k|\cdot\rho_i(c_k). \end{align*}

\noindent So, in order to calculate $p_j(i)$ for $\mc{X}(G,H\backslash G)$,
we need to know
\begin{itemize}
\item the conjugacy classes and group characters of $G$,
\item the set of double cosets of $H$ in $G$,
\item the size of the intersection of each conjugacy class and each
double coset, and
\item the decomposition of $1_H^G$ into irreducible characters.
\end{itemize}
\end{const}

This procedure has been employed when W. Kwok (1991) \cite{Kw91}
calculated the character table of $\mc{X}(O_3(q), O_2^+(q)\backslash
O_3(q))$, the scheme obtained from the action of the general
orthogonal group $O_3(q)$ acting on the sets of hyperplanes {(for
odd $q$; see \cite{Ta02} for even $q$)}. This character table
controls the character tables of association schemes coming from the
orthogonal groups on the sets of hyperplanes in the corresponding
orthogonal geometries. We remark that the character tables of the
association schemes coming from the action of $O_{2m}^{+}(q)$ on the
set of non-isotropic points are obtained by modifying the character
tables of the group $PSL(2,q)$ in exactly the same way as that of
$\mc{X}(\mc{M}^*)$ is obtained from that of $\mc{X}(PSL(2,q))$. The
reason for this may be explained as follows. Let $V$ be a
$2m$-dimensional vector space over $GF(q)$. Let $X$ be the set of
non-isotropic points corresponding to a non-singular quadratic form
of Witt index $m$. Then $|X|=q^{m-1}(q^m-1)$. The group
$G=O_{2m}^{+}(q)$ acts transitively on $X$ if $q$ is even, and it
acts transitively on each half of $X$ if $q$ is odd. Any of these
transitive permutation groups gives a symmetric association scheme
of class $q$ if $q$ is even, and class $(q+1)/2$ if $q$ is odd. It
is shown that when $m=4$ this permutation group $O_8^{+}(q)$ on $X$
(or the half of $X$) is isomorphic to the permutation group
$Gr(\mc{M}^*(q))$ on $\mc{M}^*(q)$.\vhs

The following table summarizes the results from \cite{BSS90, BSS91,
BS90, Kw91, Ta01}. In every case, the character tables of
association schemes corresponding to permutation groups $G$ on
corresponding geometries are controlled by the character table of a
`canonical' one. This happens to almost all cases (cf.
\cite{BKwS90}).\vs

\noindent\underline{Table 1. Examples of schemes from source group
1.}\vhs

\begin{center}
\renewcommand{\arraystretch}{1.2}
\begin{tabular}{|l|l|l|}
\hline  Groups & Geometries & Controlled by\\
\hline
$O_{2m}^{\pm}(q)$, $q$ even  & nonisotropic points (or lines) & $PGL(2,q)$\\
$O_{2m}^{\pm}(q)$, $q$ odd & each half of nonisotropic points & $PSL(2,q)$\\
{$O_{2m+1}(q)$}, $q$ even & $\pm$-type hyperplanes & {$PGL(2,q)/D_{2(q\mp 1)}$}\\
$O_{2m+1}(q)$, {$q$ odd} & $\pm$-type nonisotropic points & $PGL(2,q)/D_{2(q\mp 1)}$\\
\hline  $U_m(q)$ & nonisotropic points & $PGL(2,q)/Z_{q+1}$\\
\hline   $Sp_{2m}(q)$ & nonisotropic lines & $PGL(2,q)/Z_{q-1}$\\
\hline
 $P\Gamma L (n,q)$ & non-incident point-hyperplane pairs & $PGL(2,q)/Z_{q-1}$\\
\hline
\end{tabular}
\end{center}
\vs

The following table includes some examples of known Gelfand pairs
for which corresponding association schemes have been investigated.
However, the character tables of the associated commutative
association schemes are not yet known for (4).\vs

\noindent\underline{Table 2. Examples of Gelfand pairs from the
source group 2.}\vhs

\begin{center}
\renewcommand{\arraystretch}{1.2}
\begin{tabular}{|c|l|l|l|}
\hline \mbox{Labels} & \mbox{Groups}  $G$ & \mbox{Subgroups} $H$  & References\\
\hline (1) & $GL(2n, q)$ & $Sp(2n, q)$ & Klyachko \cite{Kl84}; Bannai-Kawanaka-Song \cite{BKS90}\\
\hline (2) & $GL(n,q^2)$ & $GL(n,q)$ & Gow \cite{Go84}; {Henderson \cite{He01a,He01b}}\\
\hline (3) & $GL(n,q^2)$ & $GU(n,q)$ & Gow \cite{Go84}; {Henderson \cite{He01a,He01b}} \\
\hline (4) & $GL(2n,q)$ & $GL(n,q^2)$ &
Inglis-Liebeck-Saxl \cite{ILS86}; Terras \cite{Te99}; \\
& & & Bannai-Tanaka \cite{BT03}; Henderson \cite{He03}\\
\hline (5) & $GU(2n,q^2)$ & $Sp(2n,q)$ & Inglis \cite{In88}; Henderson {\cite{He03,He01c}} \\
\hline
(6) & $Sp(4,q)$ & $Sz(q)$ & Inglis \cite{In88}; {Bagchi-Sastry \cite{BS89}; Bannai-Song \cite{BS89a}} \\
\hline
\end{tabular}
\end{center}
\vs

There are many instances where the knowledge of all ingredients in
Construction \ref{const:B} does not automatically determine the character
tables. Still there are many other examples of known Gelfand pairs,
and character tables of their corresponding association schemes need
to be determined. For example, among the association schemes
$\mathcal{X}(G, H\backslash G)$ coming from the Gelfand pairs $G$
and $H$ that appeared in the Arjeh Cohen's ``Tables of Possible
Classical Distance-Transitive Groups" (found at URL:
\href{http://www.win.tue.nl/~amc/oz/dtg/classic.html}{http://www.win.tue.nl/\~{}amc/oz/dtg/classic.html}), the character
tables of the following cases need to be calculated. Here in the
table instead of an almost simple group $G$ and its maximal subgroup
$H$, socle $S$ of $G$ and $H$ are listed.\vs

\noindent\underline{Table 3. Examples of Gelfand pairs whose
character tables need to be determined.}\vhs

\begin{center}
\renewcommand{\arraystretch}{1.2}
\begin{tabular}{|c|l|l|} \hline
$S$ & \multicolumn{1}{|c|}{type $H$} & \multicolumn{1}{|c|}{Comment} \\
\hline\hline
$PSU_6(q)$ & $SU_3(q)\times SU_3(q)$ & \cite{He03}: $1_H^G$ decomposed \\
\hline
$PSU_4(q)$ & $SU_2(q)\times SU_2(q)$ & \cite{He03}: $1_H^G$ decomposed \\
\hline
$PSU_4(q)$ & $O_4^-(q)$, $q$: odd & \cite{He03}: $1_H^G$ decomposed \\
\hline
$P\Omega_{2n}^+(q)$ & stabilizer of an $O_2^-(q)$ space & {\cite{Fu03,Fu06}}: double cosets described \\
\hline
$P\Omega_{2n}^-(q)$ & stabilizer of an $O_2^-(q)$ space & {\cite{Fu03,Fu06}}: double cosets described \\
\hline
$P\Omega_n(q)$ & stabilizer of an $O_2^-(q)$ space, $nq$: odd & ? \\
\hline
\end{tabular}
\end{center}
\vhs

Using the classification of finite simple groups, multiplicity-free
maximal subgroups of almost simple groups are getting well
understood. This will eventually lead to the complete list of
association schemes of this type which we are looking for.

\section{Character tables of $\mathcal{X}(G_2(q),O_6^{\epsilon}(q))$}\label{sec:G}
This is an example that we determine the
character table by using fission relations together with
orthogonality conditions of the character table.\vhs

Let $q$ be odd, and let $G=G_2(q)$ be the Chevalley group of type
$G_2$. Let $\Omega_1$ and $\Omega_2$ denote the sets of hyperplanes
of type $O_6^+(q)$ and $O_6^-(q)$ in the $7$-dimensional orthogonal
space over $\mbb{F}_q$. $G$ acts transitively on $\Omega _1$ and
$\Omega_2$. Let $H_1$ and $H_2$ be the one-point stabilizers of $G$
on $\Omega_1$ and $\Omega_2$, respectively. Then $H_1\simeq
SL_3(q).2$ and $H_2\simeq SU_3(q).2$. The corresponding ranks of the
permutation group are $\frac12(q+5)$ and $\frac12(q+3)$,
respectively. It is shown that $\mathcal{X}(G_2(q), \Omega_2) \simeq
\mathcal{X}(O_7(q), \Omega_2)$. The character table of
$\mathcal{X}(O_7(q), \Omega_2)$ has been constructed in
\cite{BSS90}. However, the character table of $\mathcal{X} (G_2(q),\
H_1\backslash G_2(q))$ is not isomorphic to
$\mathcal{X}(O_7(q),\,\Omega_1)$ but to its fission table. We note
that $|G|=q^6(q^2-1)(q^6-1),$ and
$|\Omega_1|=[G:H_1]=\frac12q^3(q^3+1)$ with rank($G,
\Omega_1)=\frac12(q+5)=1+$ rank$(O_7(q), \Omega_1)$.\vhs

The character table of $\mathcal{X} (G_2(q),\ SL_3(q).2\backslash G_2(q))$ and
that of $\mathcal{X}(O_7(q),\,\Omega_1)$ are, respectively, given as
follows\small
\[ \left [ \begin{array}{lllllll}
1 &2(q^3-1) &(q^2-1)(q^3-1) &q^2(q^3-1) &..... &q^2(q^3-1) &\frac12q^2(q^3-1)\\
 1 &-q^2+q-2  &q^3-q^2-q+1 &-2q^2 &..... &-2q^2 &-q^2\\
1 &2q^2-2q-2 &q^3-4q^2+2q+1 &-2q^2 &.....&-2q^2 &-q^2\\
1 & -2 & -q^2+1 & & & &\\
\vdots & \vdots & \vdots & &(q^2\chi_{ij}) &_{\substack{1\le i\le\frac12(q-1) \\ 1\le j\le\frac12(q-1)}} & \\[-2mm]
1 &-2 &-q^2+1 & & & &\\
\end{array} \right ]\]

\[\left [\begin{array}{llllll}
1 &(q^2+1)(q^3-1) &q^2(q^3-1) &\dots &q^2(q^3-1) &\frac12 q^2(q^3-1)\\
\vspace{2\jot}
1 &q^3-2q^2-1 &-2q^2 &\dots &-2q^2 &-q^2\\
\vspace{2\jot}
1 &-q^2-1 & & & &\\
\vspace{2\jot}
\vdots &\vdots & &(q^2\chi_{ij}) &_{\substack{1\le i\le\frac12(q-1) \\ 1\le j\le\frac12(q-1)}} & \\[-2mm]
\vspace{2\jot} 1 &-q^2-1 & & & &
\end{array} \right ]\]

\noindent where $\chi_{ij} \in\mathbb{Q}(\theta)\cup\mathbb{Q}(\rho),\
\theta$ and $\rho$ are the $(q+1)$th- and $(q-1)$th- root of unity,
are the entries of the character table of $\mathcal{X}
(O_3(q),\,\Omega_1)$ described by Kwok \cite{Kw91} and can be
calculated from the result in \cite[\S 6]{BSS90}.\vhs

\end{document}